\documentstyle[12pt]{article}

  \newcommand{\const}{\rm const}
  \newcommand{\Var}{\rm Var}
  
  \newcommand{\Law}{\rm Law}

\begin{document}

   \begin{center}

{\bf  Moments and tail  reciprocal  connections for  the random variables }\\

\vspace{3mm}

{\bf  having generalized  Gamma - Weibull  distributions. }\\

\vspace{4mm}

{\bf M.R.Formica, \ E.Ostrovsky, \ L.Sirota. } \\

 \end{center}

 \vspace{5mm}

\ Universit\`{a} degli Studi di Napoli Parthenope, via Generale Parisi 13, Palazzo Pacanowsky, 80132,
Napoli, Italy. \\
e-mail: \ mara.formica@uniparthenope.it \\

 \ Department of Mathematics and Statistics, Bar-Ilan University,\\
59200, Ramat Gan, Israel. \\
e-mail: \ eugostrovsky@list.ru\\

 \ Department of Mathematics and Statistics, Bar-Ilan University,\\
59200, Ramat Gan, Israel. \\
e-mail: \ sirota3@bezeqint.net \\

\vspace{5mm}

\begin{center}

{\bf Abstract.} \par

\vspace{3mm}

\end{center}

 \ We  establish the one - to one bilateral interrelations between an asymptotic behavior
 for the tail of distributions for random variables and its great moments evaluation. \par
  \ Our results generalize the famous Richter's ones. \par

\vspace{5mm}

 \par

\vspace{5mm}

\section{Introduction. Notations. Statement of problem.}

\vspace{5mm}

 \hspace{3mm} Let $ \ (\Omega = \{\omega\}, \cal{B}, {\bf P})  \ $ be certain non - trivial probability space with Expectation
 $  \ {\bf E} \ $ and Variance $ \Var. $ \par

\vspace{4mm}

 {\bf Definition 1.1.} We will say that the non - negative numerical valued random variable  (r.v.) $ \ \xi \ $  has a generalized
 Gamma - Weibull  distribution, write  $ \ \Law(\xi) \in G = G(\Gamma, W; \ Q, \ \theta, \ \gamma, \ C, \ t_0),  \ $ iff

\begin{equation} \label{main def}
T[\xi](t) \le t^{\theta} \ \exp(- C \ t^{\gamma}) \ Q(t), \ t \ge t_0 = \const > 0.
\end{equation}
 \ Here and further $ \ T[\xi](t)  = {\bf P}(\xi \ge t) \ $ is the so - called {\it tail function} for the r.v.
 $ \ \xi; \ \theta= \const > -1, \ \gamma,  C, \ t_0 = \const > 0,   \ Q = Q(t) \ $  is positive continuous function such that

$$
\forall \alpha > 0 \ \Rightarrow \lim_{t \to \infty} Q(t)/t^{\alpha} = 0.
$$

 \ To be more precisely,

\begin{equation} \label{main def prec}
T[\xi](t) \le \min(1, \ t^{\theta} \ \exp(- C \ t^{\gamma}) \ Q(t) \ ), \ t \ge  0.
\end{equation}

 \vspace{4mm}

  \ The ordinary $ \ p - th, \ p \ge 1 \ $ moment for the (non - negative) random variable $ \ \xi \ $ will be denoted by $ \ m_p = m_p[\xi]: \ $

$$
m_p [\xi] \stackrel{def}{=} {\bf E} \ \xi^p, \ p \ge 1.
$$

 \ The classical Lebesgue - Riesz $ \ ||\xi||_p \ $ norm for the r.v. $ \ \xi \ $ is defined as ordinary

$$
||\xi||_p \stackrel{def}{=} \left[ \ m_p[\xi] \ \right]^{1/p} = \left[ \ {\bf E} |\xi|^p \ \right]^{1/p}, \ p \ge 1.
$$
  \ As usually, $ \ L_p = L_p(\Omega) := \{\xi: \ ||\xi||_p < \infty \}. \ $ \par

\vspace{4mm}

 \ {\bf  Our aim in this report is to establish the bilateral exact connections between the tail behavior at infinity of the tail function
 for these variables and its moment function behavior, also at infinity.} \par

 \vspace{4mm}

\ These type estimates was applied in particular case in the theory of large deviations in the probability theory by W.Richter in articles
\cite{Richter 1}, \cite{Richter 2}, \cite{Richter 3}.  Some interest applications of these estimates in the seismology may be found in the
article of Gutenberg B.,  Richter C.F. \cite{Gutenberg Richter}. \par

\vspace{3mm}

 \ Obtained in this report estimates   complement and clarify ones  (for these variables)  in the works
\cite{Buldygin}, \cite{Capone1}, \cite{Capone2}, \cite{Ermakov}, \cite{Fiorenza2}, \cite{Fiorenza-Formica-Gogatishvili-DEA2018},
\cite{fioforgogakoparakoNAtoappear}, \cite{fioformicarakodie2017}, \cite{KosOs}, \cite{KozOsSir2017}, \cite{KosOs equivalence},
\cite{Liflyand}, \cite{Ostrovsky1}, chapters 1,2.

\vspace{5mm}

\section{ Main result.  Direct estimation.}

\vspace{5mm}

 \hspace{3mm} Suppose at first for instance for the non - negative r.v. $ \ \xi \ $  that

\begin{equation} \label{at first 1}
T[\xi](t) \le \min \left[ \ t^{\beta} \ \exp(- t), \ 1 \ \right],  \ t \ge t_0 = \const > 0, \ \beta = \const > - 1.
\end{equation}
\ We deduce using the  {\it key relation}

\begin{equation} \label{key relation}
{\bf E} |\xi|^p = p \int_0^{\infty} t^{p-1} \ T[ \xi](t) \ dt:
\end{equation}

\begin{equation} \label{template estim}
{\bf E} |\xi|^p \le t_0^p + p \ \int_{t_0}^{\infty} t^{p + \beta - 1} \ e^{-t} \ dt, \ \exists t_0 = \const > 0,  \ p \ge 1;
\end{equation}
 following

\begin{equation} \label{p estimat}
||\xi||_p \le \psi_{\beta}(p) \stackrel{def}{=} p^{1/p} \ \Gamma^{1/p}(p + \beta), \ p \ge 1.
\end{equation}

\vspace{3mm}

 \ The last relation may be rewritten as follows. Introduce the so - called Grand Lebesgue Space  $ \ G\psi_{\beta} \ $ consisting on
 all the random variables $ \ \{\zeta \} \ $ having a finite norm

 \begin{equation} \label{Gpsi norm}
 ||\zeta||G\psi_{\beta} \stackrel{def}{=} \sup_{p \ge 1} \left\{ \ \frac{||\zeta||_p}{\psi_{\beta}(p)} \ \right\} < \infty;
 \end{equation}
 then  $ \ ||\xi||G\psi_{\beta} \le 1. \ $ \par

  \vspace{3mm}

  \ The theory of these spaces is explained in many works, see e.g. \cite{Buldygin}, \cite{Capone1},  \cite{Capone2},
\cite{Ermakov}, \cite{Fiorenza2}, \cite{Fiorenza-Formica-Gogatishvili-DEA2018},
\cite{fioforgogakoparakoNAtoappear}, \cite{fioformicarakodie2017}, \cite{KosOs},
\cite{KozOsSir2017}, \cite{KosOs equivalence}, \cite{Liflyand}, \cite{Ostrovsky1}. We will use further some methods
offered in these articles. \par

\vspace{4mm}

 \hspace{3mm} Let us return to the relation (\ref{at first 1}), assuming that  it is given. \par

 \vspace{3mm}

 \  We have by virtue of the {\it key relation} \ (\ref{key relation})

\begin{equation} \label{simple mom estim}
{\bf E} \xi^p \le  p \ \int_0^{\infty} \ t^{\beta + p - 1} \ e^{-t} \ dt = p \ \Gamma(p + \beta),
\end{equation}
where as ordinary $ \ \Gamma(\cdot) \ $ denotes an ordinary Euler's Gamma function. \par

\vspace{4mm}

 \ Inversely, let the inequality (\ref{simple mom estim}) be given. We intent to evaluate the tail function for the r.v.
  $ \ \xi. \ $ \par

\vspace{3mm}

 \ One can use the classical Markov - Tchebychev's inequality for all the values of the parameters $ \ p,t \ $ such that
 $ \ p,t \ge \beta + 1 \ $

$$
T[\xi](t) \le t^{-p} \ p \ \Gamma(p + \beta).
$$
 \ One can deduce choosing the value $ \ p := t, \ $  which is asymptotically   as $ \ t \to \infty \ $ optimal:

\begin{equation} \label{tail via mom simple}
T[\xi](t) \le t^{- t} \ t \ \Gamma(t + \beta )  \le t^{-t} \ \Gamma(t + \beta + 1).
\end{equation}

\vspace{3mm}

\ We obtain  applying the  famous Stirling's formula for the values $ \ p \in(\beta + 1, \infty) \ $  and denoting

$$
c_2 := e^{1/12} \ \sqrt{2 \pi}:
$$

\vspace{3mm}

\begin{equation} \label{first tail estim}
T[\xi](t) \le c_2 \ (t + \beta)^{\beta + 1/2} \ e^{-t}.
\end{equation}

\vspace{4mm}

 \ {\bf Remark  2.1.} Note that there is a "gap" of a size  $ \  \sim t^{1/2}  \ $  as $ \ t \to \infty \ $  between the  tail estimate
(\ref{at first 1}) and  moment one (\ref{simple mom estim})
 But the estimate (\ref{first tail estim}) is asymptotically non - improvable;
for instance, for the r.v. $ \ \eta \ $ having the following tail behavior

$$
T[\eta] (t) = C_0 \ t^{\beta} \ \exp(- t), \ t > t_1 = \const > 0.
$$

\vspace{5mm}

 \hspace{3mm} Let's move on to  the rigorous considerations. Suppose that our non - negative r.v. $ \ \xi \ $ is such that there exists
  a constant $ \ \beta > - 1 \ $ and a positive  numerical valued function $ \ L = L(s), \ s > 0  \ $  such that

 \begin{equation} \label{mom condit}
 {\bf E} \xi^p \le (p + \beta)^{p + \beta} \ e^{-(p + \beta)} \ L(p + \beta).
 \end{equation}

 \  We apply  once more the Markov - Tchebychev's  inequality

$$
T[\xi](t) \le \exp [ \ - p \ln t + (p + \beta) \ln(p + \beta) - (p + \beta) + \ln L(p + \beta) \ ].
$$

 \ One can choose in the last inequality the value $ \ p = t - \beta, \ t \ge  \beta + 1. \ $ \par

 \vspace{3mm}

 \ To summarize: \par

 \vspace{4mm}

 \ {\bf Proposition  2.1.}  \ We deduce under formulated conditions,   namely, positivity of the r.v. $ \ \xi \ $  and, especially,
 under  the condition (\ref{mom condit}):

 \vspace{3mm}

\begin{equation} \label{tail estim}
T[\xi](t) \le t^{\beta}\ e^{-t} \ L(t), \ t \ge  \beta + 1.
\end{equation}

\vspace{4mm}

 \ {\bf Remark 2.2.} It is no hard to deduce the {\it lower bound} for the moment function $ \ m_p[\xi] \ $ from the {\it lower tail}
 estimate. Namely, assume now that for the non - negative random variable $ \ \xi \ $ there holds the following tail estimate

\begin{equation} \label{at first}
T[\xi](t) \ge  \ \min \left[ \ t^{\beta} \ \exp(- t), \ 1 \ \right],   \ t \ge  t_0 > 0, \ \beta = \const > - 1.
\end{equation}
\ We deduce using again the  {\it key relation}  (\ref{key relation}) and the saddle - point method  that for all the greatest values
$ \ p \ge p_0 = \const > 1 \ $

\vspace{4mm}

\begin{equation} \label{template  lower estim}
{\bf E} |\xi|^p \ge C(t_0,\beta) \ p \ \Gamma( p+ \beta), \ \exists C(t_0) \in (0,\infty).
\end{equation}

\vspace{5mm}

\section{ Main result.  Inverse  estimation.}

\vspace{5mm}

 \hspace{4mm} Inversely, let the estimate (\ref{tail estim}) be given. We want to evaluate for this r.v.  $ \ \xi \ $ its the moments function
  $ \ m[\xi]_p = {\bf E} \xi^p \ $ for all the sufficiently greatest values $ \ p, \ $ of course, under additional
  restrictions on the function $ \ L = L(t). \ $ \par
 \ We have denoting $ \ M(t) := \ln L(t), \ t \ge e;  \ z := p + \beta - 1 \ $ using again the  {\it key relation} \ (\ref{key relation})
that
$$
{\bf E} \xi^p \le p \int_0^{\infty} t^{p +\beta - 1} \ e^{-t} \ L(t) \ dt  =
$$

$$
p \int_0^{\infty} \exp ( \ z \ \ln t  - t + M(t)   \ ) \ dt =
$$

$$
p z \int_0^{\infty} \exp(z \ln z + z \ln v - zv + M(zv)) \ dv.
$$

\vspace{3mm}

 \hspace{3mm} {\it  Let us suppose now that the function   } $ \  M(\cdot) \ $ {\it is slowly varying  at infinity:}

\vspace{3mm}

\begin{equation} \label{slowly var}
\forall v > 0 \ \Rightarrow \lim_{z \to \infty} \left\{ \ \frac{M(zv)}{M(z)}\ \right\} = 1.
\end{equation}

 \vspace{4mm}

\ {\bf Proposition 3.1.} \   We deduce under this (and   previous) conditions:

$$
{\bf E} \xi^p \le  C \ p \ z \ e^{z \ln z} e^{M(z)} \int_0^{\infty} v^z \ e^{-zv} \ dv =
$$

$$
 C \ p \ \Gamma(p + \beta) \ L(p + \beta -1).
$$

\vspace{4mm}

 \ {\bf Remark 3.1.} Ii is sufficient to suppose instead the condition (\ref{slowly var}) the following one

\vspace{3mm}

\begin{equation} \label{slowly bounded}
 \sup_{z \ge e} \sup_{v \ge 1} \left\{ \ \frac{M(zv)}{M(z)}\ \right\} < \infty.
\end{equation}

\vspace{5mm}

\section{Generalizations.}

\vspace{5mm}

  \hspace{3mm} It is no hard to generalize obtained result on the case when

\begin{equation} \label{MG case}
T[\xi](t) \le t^{\theta} \ \exp \left( \ - t^{\gamma} \ \right) \ Q(t), \ \xi \ge 0,  \ t \ge t_0 = \const > 0,
\end{equation}
 \ where $ \ Q = Q(t), \ t > 0  \ $ is some positive a.e. continuous function,

$$
 \theta =\const \ge 0,\ \gamma = \const > 0.
$$

\vspace{3mm}

 \ Indeed, this case may be reduced to  the investigated above by means of the changing random variable $ \  \eta = \xi^{\gamma}, \ $ so that

$$
T[\eta](t) = {\bf P}(\xi > t^{1/\gamma}) \le t^{\theta/\gamma} \ e^{-t} \ Q \left(t^{1/\gamma} \right).
$$

\vspace{4mm}

\ {\bf Proposition 4.1.} \par

\vspace{3mm}

 \ Thus, if the function $ \ S(t) = \exp Q(t^{1/\gamma}) \ $ is positive continuous and slowly varying at infinity, then as $ \ q \in  (e, \infty) \ $

\begin{equation} \label{general moment est}
 {\bf E} \ \xi^q \le \ \frac{q}{\gamma} \ \Gamma \left( \frac{q + \theta}{\gamma}  \right) \ Q \left( \frac{q + \theta}{\gamma} \ \right).
\end{equation}

\vspace{3mm}

 \ The inverse proposition is also true in the following sense. \par

 \vspace{4mm}

 \ {\bf Proposition 4.2.}

 \vspace{3mm}

 \hspace{3mm}  Let the estimate (\ref{general moment est}) be given under at the same restrictions. Then it follows once more from 
 the Tchebychev - Markov inequality 

\vspace{3mm}

\begin{equation} \label{general tail estim}
T[\xi](t^{1/\gamma}) \le t \cdot t^{-t} \cdot \Gamma( t + \theta/\gamma) \cdot Q(t + \theta/\gamma), \ t \ge 1,
\end{equation}

\vspace{3mm}

\ and by virtue of the Stirling's approximation 

\vspace{3mm}

\begin{equation} \label{approx tail estim}
T[\xi](t^{1/\gamma}) \le c_2 \cdot t^{1/2} \cdot e^{-t} \cdot (t + \theta/\gamma)^{\theta/\gamma} \cdot Q(t + \theta/\gamma),  \ t \ge 1.
\end{equation}

\vspace{5mm}

\section{Tauberian theorems.}

\vspace{5mm}

 \hspace{3mm} Tauberian theorems are called as ordinary propositions connecting behavior of certain function, for instance, tail one,
 with some  corresponding behavior its transform,  for example, power series, Laplace, Dunkle, Fourier etc. transforms, \cite{Tauber}.
We intent to ground in this section the {\it fine} interrelations between tail behavior $ \ T[\xi](t) \ $ for the random variable $ \ \xi \ $
as $ \ t \to \infty \ $ and asymptotic behavior its moment function $ \ m_p[\xi] \ $  also as $ \ p \to \infty. \ $ \par
 \ Analogous problem for the r.v. satisfying the Kramer's condition is investigated in \cite{Ostrovsky1}, chapter 1, section 1.4., pp. \ 33 - 35.
 Some applications in the reliability theory see in  \cite{Bagdasarov}. \par

\vspace{4mm}

 \ {\bf  A. \ Direct assertion.} \par

\vspace{3mm}

 \ {\bf  Proposition 5.1.  }  Suppose that for the non - negative  r.v. $ \ \xi \ $

\begin{equation} \label{tail asympt}
\lim_{t \to \infty} \left[ \ \frac{|\ln T_{\xi}(t)|}{t} \ \right] = 1.
\end{equation}
  \ Then

\begin{equation} \label{moment asympt}
\lim_{p \to \infty} \left\{ \ \frac{||\xi||_p}{p/e}  \ \right\} = 1.
\end{equation}

 \vspace{3mm}

  \ {\bf Proof.} Let $ \ \delta \in (0,1)  \ $ be  an arbitrary fixed number. There exists a value $ \ t_0 = t_0(\delta)  \in (0,\infty) \ $
such that

$$
t \ge t_0 \ \Rightarrow T[\xi](t) \le \exp(-t(1 - \delta)).
$$

 \ We estimate for the values $ \ p > p_0 = \const > 1: \ $

$$
{\bf E} \xi^p = p \ \int_0^{\infty}  \ t^{p-1} \ T[\xi](t) \ dt \le p \int_0^{t_0} t^{p-1} \ dt +
$$

$$
p \ \int_{t_0}^{\infty} t^{p-1} \ \exp(-t(1 - \delta)) \ dt \le t_0^p + p \int_0^{\infty} t^{p-1} \ \exp( \ - t(1 - \delta) \ ) \ dt =
$$

$$
t_0^p + \frac{\Gamma(p + 1)}{(1 - \delta)^p}.
$$

 \ It follows again from the Stirling's formula that

$$
\overline{\lim}_{p \to \infty} \frac{||\xi||_p}{p/e} \le 1,
$$
and quite analogously

$$
\underline{\lim}_{p \to \infty} \frac{||\xi||_p}{p/e} \ge 1.
$$
 \ Thus, the proposition 5.1 is proved. \par

 \vspace{4mm}

 \ {\bf  B. \ Inverse proposition.} \par

\vspace{3mm}

 \ {\bf  Proposition 5.2.  }  Suppose that for the non - negative  r.v. $ \ \xi \ $

\begin{equation} \label{moment asympt inver}
\lim_{p \to \infty} \left\{ \ \frac{||\xi||_p}{p/e}  \ \right\} = 1.
\end{equation}

\ Then

\begin{equation} \label{tail inver asympt}
\lim_{t \to \infty} \left[ \ \frac{|\ln T_{\xi}(t)|}{t} \ \right] = 1.
\end{equation}

 \vspace{3mm}

  \ {\bf Proof.} {\it Upper estimate.} Let  again $ \ \delta \in (0,1)  \ $ be an arbitrary fixed number. There exists a value
 $ \ p_0  = p_0(\delta) > 1 \ $ such that

$$
\forall p \ge p_0 \ \Rightarrow \ \frac{||\xi||_p}{p/e} \le 1 + \delta,
$$
following

$$
\forall p \ge p_0 \ \Rightarrow \ ||\xi||_p   \le (1 + \delta) \cdot (p/e).
$$

 \ We apply once again the famous Tchebychev - Markov's inequality for all the sufficiently great values $ \ p \ $

 $$
 T[\xi](t) \le \frac{(1 + \delta)^p \ p^p}{e^p \ t^p} = \exp \left( \ - p \ln t + p \ln(1 + \delta) + p \ln p - p \ \right).
 $$
 \ One can choose as above $ \ p := t, \ t \ge e: \ $

$$
\ln T[\xi](t) \le - t + t \ln(1 + \delta).
$$

\ Therefore,

\begin{equation} \label{tail upper est}
\overline{\lim}_{t \to \infty} \left[ \ \frac{|\ln T_{\xi}(t)|}{t} \ \right] \le 1.
\end{equation}

\vspace{3mm}

 \ Let us deduce now the {\it lover estimate.} Given:

$$
(1- \delta)^p \ \frac{p^p}{e^p} \le {\bf E} \xi^p  \le (1+ \delta)^p \ \frac{p^p}{e^p}, \ p \ge p_0 = p_0(\delta) > e.
$$

\ Put $ \  \xi = e^{\tau}; \ $  recall that $ \ \xi > 0. \ $ Then we have for the value $ \ \lambda \ge e \ $

$$
(1- \delta)^{\lambda} \ \exp(\lambda \ \ln \lambda - \lambda) \le  {\bf E} e^{\lambda \ \tau} \le (1+ \delta)^{\lambda} \ \exp(\lambda \ \ln \lambda - \lambda).
$$

 \ The announced lower estimate

\begin{equation} \label{tail lower est}
\underline{\lim}_{t \to \infty} \left[ \ \frac{|\ln T_{\xi}(t)|}{t} \ \right] \ge 1
\end{equation}
is  grounded, up to changing variables, in particular in the monograph \cite{Ostrovsky1}, chapter 1, section 1.4; see also \cite{Bagdasarov}.\par

\vspace{5mm}

\section{Concluding remarks.}

\vspace{5mm}

 \hspace{3mm} It is interest in our opinion to deduce some multidimensional version on these result.\par

\vspace{6mm}

\vspace{0.5cm} \emph{Acknowledgement.} {\footnotesize The first
author has been partially supported by the Gruppo Nazionale per
l'Analisi Matematica, la Probabilit\`a e le loro Applicazioni
(GNAMPA) of the Istituto Nazionale di Alta Matematica (INdAM) and by
Universit\`a degli Studi di Napoli Parthenope through the project
\lq\lq sostegno alla Ricerca individuale\rq\rq .\par

\end{document}